\documentclass{amsart}

\usepackage{amsmath,amsthm,amssymb,amsfonts}
\usepackage{hyperref}

\input xy
\xyoption{all}
\usepackage{enumerate}

\makeatletter
\newcommand{\labitem}[2]{%
\def\@itemlabel{\textbf{#1}}
\item
\def\@currentlabel{#1}\label{#2}}
\makeatother

\theoremstyle{plain}
\newtheorem{thm}{Theorem}
\newtheorem{prp}[thm]{Proposition}
\newtheorem{lem}[thm]{Lemma}
\theoremstyle{remark}
\newtheorem{rmk}[thm]{Remark}

\newcommand{\Hom}{\operatorname{Hom}}

\newcommand{\Aut}{\operatorname{Aut}}

\newcommand{\Out}{\operatorname{Out}}
\newcommand{\Iso}{\operatorname{Iso}}

\newcommand{\Id}{\mathrm{Id}}

\newcommand{\res}{\operatorname{r}}

\newcommand{\tr}{\operatorname{t}}

\newcommand{\iso}{\operatorname{iso}}

\newcommand{\Mack}{\operatorname{\mathsf{Mack}}}	
\newcommand{\coMack}{\operatorname{\mathsf{coMack}}}

\newcommand{\ol}{\overline}

\makeatletter
\@tempcnta=\@ne
\def\@nameedef#1{\expandafter\edef\csname #1\endcsname}
\loop\ifnum\@tempcnta<27
  \@nameedef{C\@Alph\@tempcnta}{\noexpand\mathcal{\@Alph\@tempcnta}}
  \advance\@tempcnta\@ne
\repeat

\makeatletter
\@tempcnta=\@ne
\def\@nameedef#1{\expandafter\edef\csname #1\endcsname}
\loop\ifnum\@tempcnta<27
  \@nameedef{B\@Alph\@tempcnta}{\noexpand\mathbb{\@Alph\@tempcnta}}
  \advance\@tempcnta\@ne
\repeat

\makeatletter
\@tempcnta=\@ne
\def\@nameedef#1{\expandafter\edef\csname #1\endcsname}
\loop\ifnum\@tempcnta<27
  \@nameedef{F\@Alph\@tempcnta}{\noexpand\mathsf{\@Alph\@tempcnta}}
  \advance\@tempcnta\@ne
\repeat

\newcommand{\fgmod}[1]{#1\text{-}\mathsf{mod}}

\newcommand{\ev}{\mathrm{ev}}

\begin{document}
\title{Mislin's theorem for fusion systems  via Mackey functors}
\author{Sejong Park}
\address{School of Mathematics, Statistics and Applied Mathematics, National University of Ireland, Galway, Ireland}
\email{sejong.park@nuigalway.ie}

\date{\today}
\begin{abstract}
We state and prove a fusion system version of Mislin's theorem \cite{Mislin1990} on cohomology and control of fusion, following Symonds's proof \cite{Symonds2004} of Mislin's theorem using Mackey functors.
\end{abstract}
\maketitle

\section{Mislin's theorem for fusion systems} \label{S:Intro}

A celebrated theorem of Mislin \cite{Mislin1990} shows that an isomorphism on mod-$p$ cohomology ($p$ a prime) implies control of $p$-fusion among compact Lie groups, and in particular among finite groups. Cartan and Eilenberg's stable elements theorem \cite[XII.10.1]{CartanEilenbergBook} tells us that the mod-$p$ cohomology ring $H^*(G,\BF_p)$ of a finite group $G$ is isomorphic to the subring of the $G$-stable elements in the mod-$p$ cohomology $H^*(S,\BF_p)$ of a Sylow $p$-subgroup $S$ of $G$. In the language of fusion systems, this fact amounts to that $H^*(G,\BF_p)$ is determined by the fusion system $\CF_S(G)$ as a limit:
\[
	H^*(G,\BF_p) \cong \varprojlim_{\CF_S(G)} H^*(-,\BF_p).
\]
Here $\CF_S(G)$ is the category whose objects are the subgroups of $S$ and morphisms are the $G$-conjugation maps between them, and $H^*(-,\BF_p)$ is regarded as a contravariant functor from $\CF_S(G)$ to the category of $\BF_p$-vector spaces.

More generally, a saturated fusion system $\CF$ on a finite $p$-group $S$ is a category whose objects are the subgroups of $S$ and morphisms are group monomorphisms between them which satisfy some additional conditions formulated by L. Puig modeled on the category $\CF_S(G)$. See the book \cite{AKO} for precise definitions and more information on fusion systems, including historical background. Thus it is natural to define 
\[
	H^*(\CF,\BF_p)=\varprojlim_{\CF} H^*(-,\BF_p)
\]
for any saturated fusion system $\CF$ on $S$. With this definition Mislin's theorem for finite groups can be generalized to saturated fusion systems as follows.

\begin{thm}[Mislin's theorem for fusion systems] \label{T:Mislin-fusion}
Let $k$ be a field of characteristic~$p$. Let $\CF$ be a saturated fusion system on a finite $p$-group $S$ and let $\CE$ be a saturated subsystem of $\CF$ on the same $p$-group $S$. Suppose $H^*(\CF,k) = H^*(\CE,k)$. Then $\CF = \CE$.
\end{thm}

Recently Todea \cite{Todea2013Mislin} proved the above theorem when $p$ is odd following Benson, Grodal and Henke's new algebraic proof \cite{BensonGrodalHenke2013} of Mislin's theorem for finite groups for $p$ odd. We prove the above theorem for all primes $p$. Consequently Todea's results in \cite{Todea2013Mislin} on block algebras now hold for all primes $p$. Our approach follows closely Symonds's  proof \cite{Symonds2004} of Mislin's theorem for finite groups, which uses Mackey functors. Throughout this paper, $p$ is a fixed prime and all modules are finitely generated left modules.

\section{Simple cohomological Mackey functors and stable elements}

First we review the theory of Mackey functors with applications to fusion systems and group cohomology in mind. Specifically, we will focus on global Mackey functors. This has been done more thoroughly by D\'iaz and the author in \cite{DiazPark1}. As in  \cite{DiazPark1}, we will mainly use the language of biset functors in our discussion. See \cite{BoucBisetBook} for a systematic treatment of general biset functors, and see \cite{Webb1993} for a different approach.

Let $k$ be a commutative ring  with identity element. Let $\CI$ denote the category of finite groups and group monomorphisms. Consider the pair $M=(M^*,M_*)$ consisting of a contravariant functor $M^*\colon \CI \to \fgmod{k}$ and a covariant functor $M_*\colon \CI \to \fgmod{k}$ which have common value $M(G):=M^*(G)=M_*(G)$ on each finite group $G$ and such that $M^*(\varphi) M_*(\varphi)=\Id_{M(G)}$ for every group isomorphism $\varphi\colon G\to \varphi(G)$. Since every group monomorphism decomposes as a group isomorphism followed by an inclusion, such a pair is prescribed by a family of $k$-linear maps
\begin{gather*}
	\res^G_H := M^*(\iota_H^G) \colon M(G) \to M(H) \\
	\tr^G_H := M_*(\iota_H^G) \colon M(H) \to M(G) \\
	\iso(\varphi) := M_*(\varphi) = M^*(\varphi^{-1}) \colon  M(G) \to M(\varphi(G))
\end{gather*}
for all inclusions $\iota_H^G\colon H \to G$ and isomorphisms $\varphi\colon G\to \varphi(G)$ of finite groups, which satisfy some obvious compatibility conditions. We say that $M$ is a \emph{global Mackey functor} over $k$ if it satisfies the following additional conditions:
\begin{enumerate}
\item $\iso(c_x) = \Id_{M(G)}$ if $x\in G$.
\item (Mackey decomposition) For $H,K\leq G$,
\[
	\res^{G}_{H} \circ \tr^{G}_{K} = \sum_{x \in [H \backslash G/K]} \tr^{H}_{H\cap{}^{x}K} \circ \res^{{}^{x}K}_{H\cap {}^{x}K} \circ \iso(c_x|_K).
\]
\end{enumerate}
Here $c_x \colon G \to G$ is the conjugation map given by $c_x(u) = xux^{-1}$ for $u\in G$, and $[H \backslash G/K]$ denotes a set of double cosets $HxK$ in $G$. Note that the first condition ensures that the sum in the Mackey decomposition condition is well-defined. We say that the global Mackey functor $M$ is \emph{cohomological} if 
\[
	\tr^G_H \res^G_H = |G:H|\Id_{M(G)}
\]
for all finite groups $H\leq G$. We denote by $\Mack_k$ the category of global Mackey functors over $k$, and by $\coMack_k$ the full subcategory of cohomological global Mackey functors, where the morphism of global Mackey functors $M \to N$ is defined as the family of $k$-linear maps $M(G) \to N(G)$ for all finite groups $G$ satisfying the obvious commutativity conditions with the structure maps $\res^G_H$, $\tr^G_H$ and $\iso(\varphi)$. For example, the cohomology functor $H^*(-,k)$ with constant coefficients $k$ belongs to $\coMack_k$ where the maps $\res^G_H$, $\tr^G_H$ and $\iso(\varphi)$ are the usual restriction, transfer and isomorphism maps for group cohomology.

It is useful to view the global Mackey functors as linear functors on the bifree biset category as follows. For finite groups $G$ and $H$, let $B(H,G)$ denote the bifree double Burnside group, i.e., the Grothendieck group of isomorphism classes $[X]$ of $(H,G)$-bisets $X$ which are free as left $H$-sets and as right $G$-sets. (Note the nonstandard notation. Usually in the literature $B(H,G)$ is used for the full double Burnside group, without any freeness condition.) It is a free abelian group on the isomorphism classes of transitive bifree $(H,G)$-bisets which are of the form
\[
	H\times_{(\varphi,U)}G := H\times G/\sim,\quad (h\varphi(u),g)\sim(h,ug), h\in H, g\in G, u\in U,
\]
where $U\leq G$ and $\varphi\colon U\to H$ is a group monomorphism. The bifree biset category $\CD$ (over $k$) is the category whose objects are the finite groups and whose morphisms are given by $\Hom_\CD(G,H) = kB(H,G) := k\otimes_\BZ B(H,G)$. The composition of morphisms is induced by the tensor product of bisets: for a $(H,G)$-biset $X$ and a $(K,H)$-biset $Y$, the tensor product $Y \times_H X$ is defined as the set of $H$-orbits of $Y\times X$ where $H$ acts via $(y,x).h = (yh,h^{-1}x)$ for $x\in X$, $y\in Y$, $h\in H$. It is well-known that the global Mackey functors $M=(M^*,M_*)$ over $k$ can be viewed as $k$-linear functors $F\colon\CD\to\fgmod{k}$ and vice versa by the following corresponence:
\begin{align*}
	M^*(\varphi) &= F([G\times_{(\varphi^{-1},\varphi(G))}H]),\\
	M_*(\varphi) &= F([H\times_{(\varphi,G)}G]), \\
	F([H\times_{(\psi,U)}G]) &= M_*(\psi)M^*(\iota^G_{U}),
\end{align*}
for group monomorphisms $\varphi\colon G\to H$ and $\psi\colon U\to H$ with $U\leq G$.

As the category $\fgmod{k}$ is abelian, so is $\Mack_k$, where a sequence of global Mackey functors $L\to M\to N$ is exact if and only if its evaluation $L(G) \to M(G) \to N(G)$ at every finite group $G$ is exact. Thus we can talk about simple global Mackey functors, which we describe in the following. This is a standard fact which can be found for example in \cite[4.3]{BoucBisetBook} in a more general form.   

First we introduce the global Mackey functor $L_{H,V}$, for a finite group $H$ and a $kB(H,H)$-module $V$, given by
\[
	L_{H,V}(G) = kB(G,H) \otimes_{kB(H,H)} V
\]
for each finite group $G$ and such that the functor structure is induced by the tensor product of bisets. It is well-known that the functor $L_{H,-} \colon \fgmod{kB(H,H)} \to \Mack_k$ is left adjoint to the evaluation functor at $H$, $\ev_H \colon \Mack_k \to \fgmod{kB(H,H)}$, which sends $M$ to $M(H)$.

It is a standard technique in representation theory to mod out `things coming from below'. We adopt the idea and notation from \cite{BoucStancuThevenaz2013} as follows. For finite groups $G$, $H$, let 
\[
	I(G,H)= \sum_{K<H} kB(G,K) kB(K,H),
\]
and set $I_H = I(H,H)$ for simplicity. Let $k\ol{B}(G,H) = kB(G,H)/I(G,H)$. It turns out that $I_H$ is the submodule spanned by the basis elements $[H\times_{(U,\varphi)}H]$ with $U<H$ and $k\ol{B}(H,H) \cong k\Out(H)$. Thus we may view a $ k\Out(H)$-module as a $kB(H,H)$-module by inflation.

We say that a finite group $H$ is a \emph{minimal group} for a global Mackey functor $M$ if $M(H) \neq 0$ and $M(K) = 0$ for all $K<H$. In this case $I_H M(H) = 0$ and so $M(H)$ is a module for $k\ol{B}(H,H) \cong k\Out(H)$.

\begin{prp}[Simples in $\Mack_k$, {\cite[4.3]{BoucBisetBook}}] \label{T:simple and projective in Mack}
Let $k$ be a field.
\begin{enumerate}
\item For a finite group $H$ and a simple $kB(H,H)$-module $V$, the global Mackey functor $L_{H,V}$ over $k$ has a unique simple quotient functor $S_{H,V}$. We have $S_{H,V}(H) \cong V$. Moreover, if $V$ is a simple $k\Out(H)$-module, then $H$ is a minimal group for $S_{H,V}$.
\item If $S$ is a simple global Mackey functor over $k$ and $H$ is a finite group such that $S(H)\neq 0$, then $V=S(H)$ is a simple $kA(H,H)$-module and $S \cong S_{H,V}$. Moreover, if $H$ is a minimal group for $S$, then $V=S(H)$ is a simple $k\Out(H)$-module and the pair $(H,V)$ is uniquely determined by the isomorphism class of $S$ up to the obvious isomorphism of such pairs..
\end{enumerate}
\end{prp}

The assumption that $k$ is a field is not a real restriction because every simple Mackey functor is defined over a field. See \cite[4.4.4]{BoucBisetBook} for a proof and see also \cite[3.3]{DiazPark1} for a reduction of a part of the argument.

Simple objects in $\coMack_k$ are precisely simple objects in $\Mack_k$ which are cohomological. This is because every subfunctor and quotient functor of a cohomological Mackey functor is again cohomological. There is a particularly nice description of simple objects in $\coMack_k$ when $k$ has characteristic $p$. 

\begin{prp}[Simples in $\coMack_k$] \label{T:simple and projective in coMack}
Let $k$ be a field of characteristic $p$. Let $H$ be a finite group and let $V$ be a simple $k\Out(H)$-module. The simple global Mackey functor $S_{H,V}$ is cohomological if and only if $H$ is a $p$-group, and every simple cohomological global Mackey functor is of this form.
\end{prp}

This proposition can be proven using the explicit description of the simple global Mackey functors $S_{H,V}$ in \cite[2.6]{Webb1993}. It is also proven in a different way in \cite[5.4]{DiazPark1} as a consequence of a stronger statement, which we make a separate lemma for later use as follows. We modify $L_{H,V}$ and define the global Mackey functor $\ol{L}_{H,V}$, for a finite group $H$ and a $k\Out(H)$-module $V$, by
\[
	\ol{L}_{H,V}(G) = k\ol{B}(G,H) \otimes_{k\Out(H)} V.
\]
The functor $\ol{L}_{H,V}$ is a nonzero quotient of $L_{H,V}$ and hence $\ol{L}_{H,V}$ also has $S_{H,V}$ as a unique simple quotient. 

\begin{lem} \label{T:L-bar-cohomological}
Let $k$ be a field of characteristic $p$. If $Q$ is a finite $p$-group and $V$ is a $k\Out(Q)$-module, then the global Mackey functor $\ol{L}_{Q,V}$ is cohomological.
\end{lem}

Finally we discuss stable elements for cohomological global Mackey functors. Let $k$ be a commutative ring with identity element in which every integer prime to $p$ is invertible.  Let $M = (M^*,M_*) \in \coMack_k$. As for $H^*(-,k)$ in Section~\ref{S:Intro}, for any finite group $G$ with a Sylow $p$-subgroup $S$, we have
\[
	M(G) \cong \varprojlim_{\CF_S(G)} M^*.
\]
Therefore, for any saturated fusion system $\CF$ on a finite $p$-group $S$, it is natural to define
\[
	M(\CF) = \varprojlim_{\CF} M^*.
\]
When $\CF=\CF_S(G)$, the group $G$ viewed as the obvious $(S,S)$-biset determines the value $M(G)$ as $M(G) \cong [G] \cdot M(S)$. In general, every saturated fusion system $\CF$ on $S$ has an idempotent element $\omega_\CF \in \BZ_{(p)}B(S,S)$ which plays the role of the $(S,S)$-biset $G$ for $\CF_S(G)$, called the \emph{characteristic idempotent} of $\CF$. See \cite{Ragnarsson2006} for more details. We record the following well-known fact for later use.

\begin{lem}[{\cite[2.8, 2.9]{DiazSpectral}}] \label{T:ev-at-F exact}
Let $k$ be a commutative ring with identity element in which every integer prime to $p$ is invertible. Let $M=(M^*,M_*)$ be a cohomological global Mackey functor and let $\CF$ be a saturated fusion system with characteristic idempotent $\omega_\CF$.
\begin{enumerate}
\item $M(\CF) \cong \omega_\CF \cdot M(S)$.
\item The evaluation functor
\[
	\coMack_k \to \fgmod{k},\quad M\mapsto M(\CF),
\]
is exact.
\end{enumerate}
\end{lem}

\section{A proof of Mislin's theorem} \label{S:Mislin}

Now we prove Theorem~\ref{T:Mislin-fusion} following Symonds's argument in \cite{Symonds2004} closely. Recall that, in Theorem~~\ref{T:Mislin-fusion}, $k$ is a field of characteristic $p$, $\CF$ is a saturated fusion system on a finite $p$-group $S$, and $\CE$ is a saturated subsystem of $\CF$ on $S$. Symonds \cite[1.2]{Symonds2004} showed that every simple object in $\coMack_k$ occurs as a composition factor of $H^*(-,k)$. If 
\[
	0 \to L \to M \to N \to 0
\]
is a short exact sequence in $\coMack_k$, by Lemma~\ref{T:ev-at-F exact} we have a commutative diagram with exact rows:
\[
\xymatrix{
	& 0 \ar[r] & L(\CF) \ar[r]\ar@{^{(}->}[d] & M(\CF) \ar[r]\ar@{^{(}->}[d] &N(\CF) \ar[r]\ar@{^{(}->}[d] & 0 \\
	& 0 \ar[r] & L(\CE) \ar[r] & M(\CE) \ar[r] &N(\CE) \ar[r] & 0, \\
}
\]
and an easy diagram chase shows that
\[
	M(\CF) = M(\CE) \iff L(\CF) = L(\CE), N(\CF) = N(\CE).
\]
Therefore the following statements are equivalent:
\begin{enumerate}
\item $H^*(\CF,k) = H^*(\CE,k)$.
\item $M(\CF) = M(\CE)$ for every simple $M \in \coMack_k$.
\item $M(\CF) = M(\CE)$ for every $M \in \coMack_k$.
\end{enumerate}
Since $\ol{L}_{Q,k\Out(Q)} = k\ol{B}(-,Q)$ belongs to $\coMack_k$ if $Q$ is a finite $p$-group by Lemma~\ref{T:L-bar-cohomological}, the assumption of the theorem that $H^*(\CF,k) = H^*(\CE,k)$ implies that
\[
	k\ol{B}(\CF,Q) = k\ol{B}(\CE,Q)\quad\forall\, Q\leq S,
\]
which is equivalent to that
\[
	\omega_\CF \cdot k\ol{B}(S,Q) = \omega_\CE \cdot k\ol{B}(S,Q)\quad\forall\, Q\leq S
\]
by Lemma~\ref{T:ev-at-F exact}. Ragnarsson \cite[5.2]{Ragnarsson2006} gives a $k$-basis for $\omega_\CF \cdot kB(S,Q)$, and from it one can easily deduce a $k$-basis for $\omega_\CF \cdot k\ol{B}(S,Q)$. In particular, the $k$-dimension of $\omega_\CF \cdot k\ol{B}(S,Q)$ is equal to the number of  group isomorphisms $\varphi\colon Q\to L$ with $L$ a subgroup $S$, up to composition with $\CF$-isomorphisms. Since $\Aut_\CF(L)$ acts freely on the set $\Iso(Q,L)$ of isomorphisms from $Q$ to $L$ by composition and any group isomorphism $\varphi\colon Q\to L\leq S$ gives a bijection from $\Iso(Q,L)$ to $\Aut(L)$ by sending $\psi\in\Iso(Q,L)$ to $\psi\circ\varphi^{-1}$, we have
\begin{align*}
	\dim_k \omega_\CF \cdot k\ol{B}(S,Q)= \sum_{\substack{Q \cong L\leq S\\\text{up to $\CF$-iso.}}} \frac{|\Iso(Q,L)|}{|\Aut_\CF(L)|} = \sum_{\substack{Q \cong L\leq S\\\text{up to $\CF$-iso.}}} \frac{|\Aut(L)|}{|\Aut_\CF(L)|},
\end{align*}
where $L$ runs over the subgroups of $S$ which are isomorphic to $Q$, taken up to $\CF$-conjugacy. Summing up these dimensions over the subgroups $Q$ of $S$ up to isomorphism, and comparing it with that for $\CE$, we get an equality
\[
	\sum_{\substack{Q \leq S\\\text{up to iso.}}} \sum_{\substack{Q \cong L\leq S\\\text{up to $\CF$-iso.}}} \frac{|\Aut(L)|}{|\Aut_\CF(L)|} = \sum_{\substack{Q \leq S\\\text{up to iso.}}} \sum_{\substack{Q \cong L\leq S\\\text{up to $\CE$-iso.}}} \frac{|\Aut(L)|}{|\Aut_\CE(L)|}.
\]
Since the right-hand side has at least as big summands and at least as many terms as the left-hand side does, the above equality forces that $\CF$ and $\CE$ have the same number of isomorphism classes of subgroups of $S$ and the same automorphism groups $\Aut_\CF(L)=\Aut_\CE(L)$ for each subgroup $L$ of $S$. By Alperin's fusion theorem, it follows that $\CF=\CE$.

\begin{rmk}
Note that the above proof is all algebraic except for the quoted fact \cite[1.2]{Symonds2004}, which has been proven using stable homotopy theory, including (proven) Segal's Burnside ring conjecture. See \cite[\S 4]{Symonds2004} for more details. Thus an algebraic proof of \cite[1.2]{Symonds2004} would give an algebraic proof of Theorem \ref{T:Mislin-fusion}.
\end{rmk}

\end{document}